\documentclass[12pt]{article}
\usepackage[latin1]{inputenc}
\usepackage{amssymb, amsmath, amsthm, bm, epsfig}
\usepackage{graphicx, color}
\usepackage{pdflscape}
\usepackage{lscape}
\usepackage{rotating}

\newcommand {\tr} {\mbox{tr}} 
     
\newcommand{\vvec}{\textrm{vec}}
\newcommand {\J}{\widetilde{\!\!\widetilde{J}}}
\newcommand {\inic} {\hspace{0.5cm}}

\title{Improved hypothesis testing in a general multivariate elliptical model}  
\author
{
   Tatiane F. N. Melo\\
     {\it \footnotesize Institute of Mathematics and Statistics, Federal University of Goi\'as, Brazil} \vspace{-0.2cm}\\
     {\footnotesize email: {\tt tmelo@ufg.br}} \\ \\
   Silvia L. P. Ferrari\\
   {\it \footnotesize Departament of Statistics, University of S\~ao Paulo, Brazil} \vspace{-0.2cm}\\
     {\footnotesize email: {\tt silviaferrari@usp.br}} \\ \\
   Alexandre G. Patriota\\
     {\it \footnotesize Departament of Statistics, University of S\~ao Paulo, Brazil} \vspace{-0.2cm}\\
     {\footnotesize email: {\tt patriota@ime.usp.br}} \\ \\
}
\date{}

\begin{document}
\maketitle

%------------------------------------- RESUMO --------------------------------------------
{\footnotesize 
\noindent{\bf Abstract:} 

\noindent
This paper investigates improved testing inferences under a general multivariate elliptical regression model. The model is very flexible in terms of the specification of the mean vector and the dispersion matrix, and of the choice of the error distribution. The error terms are allowed to follow a multivariate distribution in the class of the elliptical distributions, which has the multivariate normal and Student-t distributions as special cases. We obtain Skovgaard's adjusted likelihood ratio statistics and Barndorff-Nielsen's adjusted signed likelihood ratio statistics and we conduct a simulation study. The simulations suggest that the proposed tests display superior finite sample behavior as compared to the standard tests. Two applications are presented in order to illustrate the methods.

\vspace{0.1cm}
}

%------------------------------------- INTRODUÇÃO --------------------------------------------
\section{Introduction}
\label{sec1}

\inic Likelihood inference is usually based on the first order asymptotic theory, which can lead to inaccurate inference when the sample is small. In general, this is the case of the signed likelihood ratio test, whose statistic has asymptotic standard normal distribution under the null hypothesis, with an error of order $n^{-1/2}$, where $n$ is the size sample. In order to improve this approach, Barndorff-Nielsen (1986) proposed a new test statistic, that under the null hypothesis is asymptotically standard normal distributed with error of order $n^{-3/2}$. Barndorff-Nielsen's adjustment is applied when the parameter of interest is scalar. Skovgaard (2001) developed an extension of this adjustment for the multidimensional case. These adjustments require a suitable ancillary statistic such that, in conjunction with the maximum likelihood estimator, they must constitute a sufficient statistic for the model. It is difficult or even impossible to find an appropriate ancillary for some statistical models (Pe\~na et al, 1992). In this paper, we obtain Barndorff--Nielsen's and Skovgaard's adjustments in a general multivariate elliptical model using an approximate ancillary statistic. We perform a simulation study that suggests that the modified tests have type I probability error closer to the nominal level the original tests in small and moderate-sized samples.

A general multivariate elliptical model was introduced by Lemonte and Patriota (2011). It considers that the mean vector and the dispersion matrix are indexed by the same vector of parameters. Multiple linear regressions, multivariate nonlinear regressions, mixed-effects models (Verbeke and Molenberghs, 2000), errors-in-variables models (Cheng and Van Ness, 1999) and log-symmetric regression models (Vanegas and Paula, 2015) are special cases of this general multivariate elliptical model.

The elliptical family of distributions includes the multivariate normal as well as many other important distributions such as the multivariate Student $t$, power exponential, contaminated normal, Pearson II, Pearson VII, and logistic distributions, with heavier or lighter tails than the multivariate normal distribution. The random vector ${\bm Y}$ ($q\times 1$) has a multivariate elliptical distribution with location parameter ${\bm \mu}$ ($q\times 1$) and a positive definite scatter matrix $\Sigma$ ($q\times q$) if its density function exists and is given by 
\begin{equation*}\label{dens}
f_{{{\bm Y}}}({{\bm y}}) = |{\Sigma}|^{-1/2} g\bigl(({{\bm y}} - {{\bm \mu}})^{\top}{\Sigma}^{-1}({{\bm y}} - {{\bm \mu}})\bigr),
\end{equation*}
where $g:[0,\infty)\to(0,\infty)$ is called the density generating function and it is such that $\int_{0}^{\infty}u^{\frac{q}{2}-1} g(u) du < \infty$. We will denote ${{\bm Y}}\sim El_{q}({{\bm \mu}},{\Sigma}, g) \equiv El_{q}({{\bm \mu}},{\Sigma})$. The characteristic function is $\psi({\bm t}) = \mbox{E}(\exp(i {\bm t}^\top {\bm Y}))=\exp(i {\bm t}^\top {\bm \mu})\varphi({\bm t}^\top \Sigma {\bm t})$, where ${\bm t} \in \mathbb{R}^q$ and $\varphi: [0,\infty) \to \mathbb{R}$. Then, if $\varphi$ is twice differentiable at zero, we have that $\mbox{E}({\bm Y}) = {\bm \mu}$ and $\mbox{Var}({\bm Y}) = \xi \Sigma$, where $\xi = \varphi'(0)$. We assume that the density generating function $g$ does not have any unknown parameter, which implies that $\xi$ is a known constant. In this case, when ${{\bm \mu}} = {\bm 0}$ and ${\Sigma} = {I}_{q}$, where ${I}_{q}$ is the $q\times q$ identity matrix, we obtain the spherical family of densities; for more details see Fang et al.~(1990).

This paper is organized as follows. Section \ref{sec2} introduces the general elliptical model. Section \ref{sec3} contains our main results, namely explicit formulas for  Barndorff-Nielsen's and Skovgaard's adjustments in the general elliptical model. Section \ref{sec4} presents a simulation study on the finite sample behavior of the standard and signed likelihood ratio tests and their modified counterparts. Our simulation results show that the unmodified tests tend to be liberal and their modified versions developed in this paper are much less size-distorted.  Section \ref{sec5} presents two real data applications. Finally, Section \ref{sec6} concludes paper. Technical details are compiled in the Appendix.

%------------------------------------- O MODELO --------------------------------------------
\section{The model}
\label{sec2}

\inic The general multivariate elliptical model defined by Lemonte and Patriota (2011) considers $n$ independent random vectors ${\bm Y}_{1}, {\bm Y}_{2},\ldots, {\bm Y}_{n}$ modeled by the following equation:
\begin{equation}\label{MainModel1} 
{\bm Y}_{i} = {\bm \mu}_{i}({\bm \theta}) + {\bm e}_{i},\quad i = 1,\ldots,n,
\end{equation}
with ${\bm e}_{i} \stackrel{ind}{\sim} El_{q_{i}}({\bm 0},{\Sigma}_i({\bm \theta}))$, where ``$ \stackrel{ind}{\sim}$'' means ``independently distributed as'', ${\bm \mu}_{i}({\bm \theta}) = {\bm \mu}_{i}$ is the location parameter and ${\Sigma}_i({\bm \theta}) = {\Sigma}_i$ is the positive definite scatter matrix\footnote{Note that ${\bm \mu}_{i}$ and ${\Sigma}_i$ may depend on covariates 
associated with the $i$th observed response ${\bm Y}_{i}$. The covariates may have components in common.}  
We can write
\begin{equation}\label{MainModel}
{\bm Y}_{i} \stackrel{ind}{\sim}El_{q_{i}}({\bm \mu}_{i}, {\Sigma}_i),\quad i = 1,\ldots,n.
\end{equation}
\noindent Both ${\bm \mu}_{i}$ and ${\Sigma}_i$ have known functional forms.
Additionally, ${\bm \theta}$ is a vector of unknown parameters, with ${\bm \theta} \in {\bm \Theta} \subseteq \mathbb{R}^p$ (where $p<n$ is fixed). 

For the ${q}$-variate normal distribution, $N_{q}({\bm{\mu}}, {\Sigma})$, the density generating function is  $g(u) = e^{-u/2}/(\sqrt{2 \pi})^{q}$. For the ${q}$-variate Student $t$ distribution with $\nu$ degrees of freedom, $t_{q}({\bm{\mu}}, {\Sigma}, \nu)$, we have $g(u) = \Gamma\left((\nu + q)/2\right) \pi^{-q/2} \nu^{-q/2} (1 + u/\nu)^{-(\nu + q)/2} / \Gamma(\nu/2) $. Additionally, for the $q$-variate power exponential, $PE_{q}({\bm{\mu}}_i, {\Sigma}, \lambda)$, with shape parameter $\lambda>0$, we have $g(u) = \lambda \Gamma(q/2) 2^{-q/(2 \lambda)} \pi^{-q/2} e^{-u^{\lambda}/2} /$ $\Gamma(q / (2 \lambda))$. 

For this general model, Lemonte and Patriota (2011) proposed diagnostic tools and Melo et al.~(2015) obtained the second-order bias of the maximum likelihood estimator and conducted some simulation studies, which indicate that the proposed bias correction is effective. 

The log-likelihood function associated with (\ref{MainModel}) is given by 
\begin{equation} \label{log-likelihood}
\ell({\bm \theta}) = \sum_{i=1}^n \ell_{i}({\bm \theta}),
\end{equation}
where $\ell_{i}({\bm \theta}) = -\frac{1}{2} \log{|{\Sigma}_i|} + \log g(u_i)$, $u_{i} = {\bm z}_i^{\top}{\Sigma}_i^{-1}{\bm z}_i$ and ${\bm z}_i = {\bm Y}_i - {\bm \mu}_i$. The dependence of $\ell({\bm \theta})$ on ${\bm \theta}$ enters through ${\bm{\mu}}_i=\mu_i({\bm \theta})$ and $\Sigma_i=\Sigma_i({\bm \theta})$. 
We assume regularity conditions for the asymptotic theory of maximum likelihood estimation and likelihood ratio tests; see Severini (2000, \S \ 3.4). The model must be identifiable and it must be guaranteed that the first four derivatives of $(1/n)\ell(\bm{\theta})$ with respect to $\bm{\theta}$ exist, are bounded by integrable functions, and converge almost surely for all $\bm{\theta}$. The conditions of the Lindeberg-Feller Theorem (or the Liapounov Theorem) must be valid for the score function to converge in distribution to a normal distribution (Sen and Singer, 1993, p.108). These conditions impose restrictions on the sequences
$\{\bm{\mu}_i\}_{i\geq 1}$ and $\{\Sigma_i\}_{i\geq 1}$ that will not be detailed here.

Maximum likelihood estimation of the parameters  can be carried out by numerically maximizing the log-likelihood function (\ref{log-likelihood}) through an iterative algorithm such as the Newton--Raphson, Fisher scoring, EM or BFGS. Our numerical results were obtained using the library function {\tt MaxBFGS} in the {\tt Ox} matrix programming language (Doornik, 2013).

%------------------------------------- AJUSTE DE SKOVGAARD --------------------------------------------
\section{Modified likelihood ratio tests}
\label{sec3}

\inic Consider the vector of unknown parameters ${\bm \theta} = ({\bm \psi}^\top, {\bm \omega}^\top)^\top\in \mathbb{R}^p$, with ${\bm \psi}\in \mathbb{R}^q$ being the vector of parameters of interest and ${\bm \omega}\in \mathbb{R}^{p-q}$ being the vector of nuisance parameters. The null and alternative hypotheses of interest are, respectively: ${\cal H}_{0}: {{\bm \psi}} = {{\bm \psi}}^{(0)}$ and ${\cal H}_{1}: {{\bm \psi}} \neq {{\bm \psi}}^{(0)}$, where  ${{\bm \psi}}^{(0)}$ is a known $q$-vector. The maximum likelihood estimator of ${\bm \theta}$ is denoted by $\bm{\widehat{\theta}} = (\bm{\widehat{\psi}}, \bm{\widehat{\omega}}^\top)^\top$  and the maximum likelihood estimator of ${\bm \theta}$ under the null hypothesis, by $\bm{\widetilde\theta} = ({\bm{\psi}}^{(0)}, \bm{\widetilde\omega} ^\top)^\top$. We use ``\ $\widehat{}$\ " and ``\ $\widetilde{}$\ " for matrices and vectors to indicate that they are computed at $\widehat{\bm{\theta}}$ and $\widetilde{\bm{\theta}}$, respectively. The standard likelihood ratio statistic for testing ${\mathcal{H}}_{0}$ against ${\mathcal{H}}_1$ is
\begin{eqnarray*}\label{E.4.11}
LR = 2\:\left\{\ell(\bm{\widehat{\theta}}) - \ell(\bm{\widetilde{\theta}})\right\}.
\end{eqnarray*}
Under regularity conditions (Severini, 2000, Section 3.4), $LR$ converges in distribution to $\chi^2_q$ when ${\cal H}_{0}: {{\bm \psi}} = {{\bm \psi}}^{(0)}$ holds.

When the parameter of interest, $\psi$, is one-dimensional, i.e., $q=1$,  the signed likelihood ratio statistic,
\begin{eqnarray*}\label{rS}
r = {\rm sgn}\left(\widehat{\psi} - {\psi}^{(0)} \:\right) \sqrt{2\:
\left(\ell(\widehat{\bm{\theta}}) - \ell(\widetilde{\bm{\theta}})\right)},
\end{eqnarray*}
may be employed. Under regularity conditions, $r$ converges in distribution to a standard normal distribution when ${\cal H}_{0}: \psi=\psi^{(0)}$ is true. In addition to two-sided tests (${\mathcal{H}}_{0}: {\psi} = {\psi}^{(0)} \:\:\: {\rm against} \:\:\: {\mathcal{H}}_{1}: {\psi} \neq {\psi}^{(0)}$), one may use the signed likelihood ratio statistic to test one-sided hypotheses such as ${\mathcal{H}}_{0}: {\psi} \geq {\psi}^{(0)} \:\:\: {\rm against} \:\:\: {\mathcal{H}}_{1}: {\psi} < {\psi}^{(0)}$ and ${\mathcal{H}}_{0}: {\psi} \leq {\psi}^{(0)} \:\:\: {\rm against} \:\:\: {\mathcal{H}}_{1}: {\psi} > {\psi}^{(0)}$. 

Barndorff-Nielsen (1986) proposes a modified version of $r$ that intends to better approximate the signed likelihood ratio statistic by the standard normal distribution. It can be difficult to obtain the modified statistic because one needs to obtain an appropriate ancillary statistic and derivatives of the log-likelihood function with respect to the data. By ``ancillary statistic" we mean a statistic, say ${\bm a}$, whose distribution does not depend on the unknown parameter $\bm{\theta}$, and such that $(\bm{\widehat{\theta}},{\bm a})$ is a minimal sufficient statistic for the model. If $(\bm{\widehat{\theta}},{\bm a})$ is sufficient, but not minimal sufficient, Barndorff-Nielsen's results still hold; see, Severini (2000, \S \ 6.5).  
Sufficiency implies that the log-likelihood function depends on the data only through $(\bm{\widehat{\theta}}, {\bm a})$, and we then write $\ell(\bm{\theta};\bm{\widehat \theta},{\bm a})$. The derivatives of  $\ell(\bm{\theta};\bm{\widehat \theta},{\bm a})$ with respect to the data and the parameter vector are
\begin{eqnarray*}\label{E.17}
\bm{\ell}'(\bm{\theta};\bm{\widehat \theta},{\bm a}) = \frac{\partial \ell(\bm{\theta};\bm{\widehat \theta},{\bm a})}
{\partial \bm{\widehat \theta}}, \ \ 
U'(\bm{\theta};\bm{\widehat \theta},{\bm a}) = \frac{\partial^2 \ell(\bm{\theta};\bm{\widehat \theta},
{\bm a})}{\partial \bm{\widehat \theta} \partial{\bm{\theta}}^\top}, \ \ {\rm and} \ \ 
J(\bm{\theta};\bm{\widehat \theta},{\bm a}) = -\frac{\partial^2 \ell(\bm{\theta};\bm{\widehat \theta},
{\bm a})}{\partial \bm{\theta} \partial{\bm{\theta}}^\top}. 
\end{eqnarray*}

The modified signed likelihood ratio statistic is
\begin{eqnarray*}\label{rstar}
r^* = r - \frac{1}{r} \log \gamma,
\end{eqnarray*}
with
\begin{equation}\label{E.19}
\gamma = |\widehat{J}\:|^{1/2} |{\widetilde U}'|^{-1} |{\widetilde J}_{\bm{\omega\omega}}|^{1/2} 
\frac{r}{[({\widehat {\bm{\ell}}}'- {\widetilde {\bm{\ell}}}')^\top ({\widetilde U}')^{-1}]_{\psi}}, 
\end{equation}
where $J = J(\bm{\theta};\bm{\widehat \theta},{\bm a})$ is the observed information matrix and $J_{\bm{\omega\omega}}$ is the lower right submatrix of $J$ corresponding to the nuisance parameter $\bm{\omega}$. Here, $[\bm{v}]_{\psi}$ denotes the element of the vector $\bm{v}$ that corresponds to the parameter of interest $\psi$. The quantities ${\widehat {\bm{\ell}}}'={\bm{\ell}}'(\bm{\widehat \theta};\bm{\widehat \theta},{\bm a})$, ${\widetilde {\bm{\ell}}}'={\bm{\ell}}'(\bm{\widetilde \theta};\bm{\widehat \theta},{\bm a})$ and ${\widetilde U}'= U'(\bm{\widetilde \theta};\bm{\widehat \theta},{\bm a})$ are computed as described above. 

Barndorff-Nielsen's $r^*$ statistic is only useful when testing a one-dimensional hypothesis. However, in practical applications it is often the case that the null hypothesis involves several parameters. As an example, we mention the test for treatment effects in linear mixed models; see Verbeke and Molenberghs (2000, \S \ 6.2). The seminal work of Skovgaard (2001) extended Barndorff-Nielsen's (1986) results to the multiparameter test situation. He proposed two modified, asymptotically equivalent, likelihood ratio statistics, which can be seen as multiparameter versions of $r^*$. The modified statistics are derived from Barndorff-Nielsen's work and share similar properties with $r^*$; see Skovgaard (2001, \S \ 5, and eq. (8)-(10)). Skovgaard's modified likelihood ratio statistics are given by

\begin{eqnarray*}\label{E.4.1}
LR^* = LR \left(1 - \frac{1}{LR} \log\rho \right)^2
\end{eqnarray*}
and
\begin{eqnarray*}\label{E.4.16}
LR^{**} = LR - 2\log\rho,
\end{eqnarray*}
with
\begin{eqnarray}\label{E.4.17}
\rho = |\widehat{J}\:|^{1/2} |{\widetilde U}'|^{-1} |
{\widetilde J}_{\bm{\omega\omega}}|^{1/2} |
\:{\J}_{\bm{\omega\omega}}|^{-1/2} |\:{\J}\:|^{1/2} 
\frac{\{{\widetilde {\bm{U}}}^{\top} {\J}^{\: -1} {\widetilde {\bm{U}}}\}^{p/2}}
{LR^{q/2 - 1} ({\widehat {\bm{\ell}}}'- {\widetilde {\bm{\ell}}}')^{\top} 
({\widetilde U}')^{-1} {\widetilde{\bm{U}}}},
\end{eqnarray}
where ${\bm{U}}$ is the score vector and $\:\:{\J} = J(\bm{\widetilde{\theta}}; \bm{\widetilde{\theta}}, {\bm a})$, and ${\:\:{\J}}_{\bm{\omega \omega}}$ is the lower-right sub-matrix of ${\:\:{\J}}$ related to the nuisance parameters $\bm{\omega}$. Although the statistic $LR^*$ is non-negative and reduces to ${r^*}^2$ when $q=1$, the second version, $LR^{**}$, seems to be numerically more stable  and is naturally attained from theoretical developments. These statistics approximately follow the asymptotic reference distribution (${\cal X}^2_q$ distribution) with high accuracy under the null hypothesis (Skovgaard, 2001). 
In fact, recent simulation studies suggest that Barndorff-Nielsen's and Skovgaard's statistics considerably improve small-sample inference; see, for example, Brazzale and Davison (2008), Lemonte and Ferrari (2011), Ferrari and Pinheiro (2014), Guolo (2012) and Cribari-Neto and Queiroz (2014).

We now turn to the general elliptical model. Let ${\bm a} = ({\bm a}_1^\top, {\bm a}_2^\top$, $\ldots, {\bm a}_n^\top)^\top$, with
\begin{equation*}\label{E.20}
{\bm a}_i = {\widehat P}_i^{-1}\left({\bm Y}_i - {\bm{\widehat \mu}_i}\right),
\end{equation*}
where $P_i\equiv P_i(\bm\theta)$ is a lower triangular matrix such that $P_i P_i^\top = \Sigma_i$ is the Cholesky decomposition of $\Sigma_i$ for all $i=1, \ldots, n$. 
From Slutsky's theorem, it follows that $\bm{a}_i$ converges in distribution to $El_{q_{i}}(\bm{0}, I_{q_i})$, since ${\widehat P}_i$ and $\widehat{\bm{\mu}}_i$ converge in probability to $P_i$ and $\bm{\mu}_i$, respectively. Additionally, it can be shown that any fixed number of the $\bm{a}_i$'s are asymptotically independent, and hence their joint asymptotic distribution is free of unknown parameters. 
Also, it follows from Neyman's Factorization Theorem that $(\widehat{\bm{\theta}},\bm{a})$ is a sufficient statistic, since the log-likelihood function can be written as 
\begin{equation}\label{E.ell}
\ell(\bm{\theta};\widehat{\bm{\theta}}, \bm{a}) = \sum_{i=1}^n \left\{ - \frac{1}{2} \log |\Sigma_i | + 
\log g[({\widehat P}_i \bm{a}_i + \bm{\widehat \mu}_i - \bm{\mu}_i)^\top \Sigma_i^{-1} ({\widehat P}_i \bm{a}_i + {\bm{\widehat \mu}_i - \bm{\mu}_i}) ] \right\},
\end{equation}
where the dependence\footnote{The dependence on $\bm{\theta}$ is omitted here and in the sequel for the sake of readability.} on $\bm{\theta}$ is through $\bm{\mu}_i$ and $\Sigma_i$.  Hence, we will use $\bm{a}$ as an approximate ancillary statistic. The use of an approximate ancillary statistic in connection with Barndorff-Nielsen's $r^{*}$ statistic can be found, for example, in Fraser, Reid and Wu (1999); see also Severini (2000, \S \ 7.5.3). Formulas (\ref{E.19}) and (\ref{E.4.17}), and the modified statistics $r^{*}$, $LR^{*}$, and $LR^{**}$ can be computed from (\ref{E.ell}). Details are presented in the Appendix.

Let
\[
{\bm d}_{i(r)}  = \frac{\partial {\bm \mu}_i}{\partial \theta_{r}}, \quad
{\bm d}_{i(sr)} = \frac{\partial^2 {\bm \mu}_i}{\partial \theta_{s} \partial\theta_{r}}, \quad 
{C}_{i(r)}      = \frac{\partial {\Sigma}_i}{\partial \theta_{r}}, \quad 
{C}_{i(sr)}     = \frac{\partial^2 {\Sigma}_i}{\partial \theta_{s} \partial\theta_{r}},
\]
and\[
{A}_{i(r)} = -{\Sigma}_i^{-1} {C}_{i(r)}{\Sigma}_i^{-1},
\]
for $r, s = 1, \ldots, p$. The score vector and the observed information matrix for ${\bm \theta}$ can be written as 
\begin{equation}\label{ScoreFisher2}
{\bm{U}} = {F}^{\top}{H}{\bm s}, \quad J = T^{\top} \Sigma^{-1} D + G,
\end{equation}
respectively, with ${F} = \left({F}_1^\top, \ldots, {F}_n^\top\right)^{\top}$, $H = \mbox{block-diag}\left\{H_1, \ldots, H_n\right\}$, ${\bm s} = ({\bm s}_1^\top, \ldots, {\bm s}_n^\top )^{\top}$, $T = \left({T}_1^\top, \ldots, {T}_n^\top\right)^{\top}$, $\Sigma^{-1} = \mbox{block-diag}\left\{ \Sigma_1^{-1},\ldots, \Sigma_n^{-1} \right\}$, ${D} = \left({D}_1^\top, \ldots, {D}_n^\top\right)^{\top}$,
wherein
\begin{equation*} \label{matrix-score}
F_i =
  \begin{pmatrix}
    {D}_i\\
    {V}_i\\
  \end{pmatrix},\quad
{H}_i =
  \begin{bmatrix}
    {\Sigma}_i & {0}\\
    {0} & 2{\Sigma}_i\otimes {\Sigma}_i
  \end{bmatrix}^{-1}, \quad
{\bm s}_i =
  \begin{bmatrix}
    v_{i}{\bm z}_i\\
    -\vvec({\Sigma}_i - v_{i} {\bm z}_i {\bm z}_i^{\top})
  \end{bmatrix},
\end{equation*}
where the ``vec" operator transforms a matrix into a vector by stacking the columns of the matrix, ${D}_i = ({\bm d}_{i(1)}, \ldots, {\bm d}_{i(p)})$, ${V}_i = (\vvec({C}_{i(1)}),\ldots, \vvec({C}_{i(p)}))$, ${T}_i = ({\bm T}_{i(1)}, \ldots, {\bm T}_{i(p)})$, $v_{i} = -2 W_{g}(u_{i})$ and $W_g(u) = \mbox{d} \log g(u)/\mbox{d} u$. The $(r,s)$-th elements of $G$ and $E_i$ are given by $G_{rs}$ and $E_{i(rs)}$, respectively. The quantities ${\bm T}_{i(r)}$, $G_{rs}$, and $E_{i(rs)}$ are given in the Appendix.  
The symbol ``$\otimes$'' indicates the Kronecker product. 

We have
\begin{eqnarray}\label{Deriv}
\widehat{\bm{\ell}}'   = \widehat{R}^{\top} \widehat{\Sigma}^{-1} \widehat{{\bm z}}^*,  \ \
\widetilde{\bm{\ell}}' = \widehat{R}^{\top} \widetilde{\Sigma}^{-1} \widetilde{{\bm z}}^*,  \ \
\widetilde{U}'    = \widetilde{Q}^{\top} \widetilde{\Sigma}^{-1} \widehat{R}, \ \
\J                = \ {\widetilde{\!\!\widetilde{{T}}}}^{\top} \widetilde{\Sigma}^{-1} \widetilde{{D}} + \ {\widetilde{\!\!\widetilde{G}}},
\end{eqnarray}
where $\widehat{R}$, $\widehat{{\bm z}}^*$, $\widetilde{{\bm z}}^*$, $\widetilde{Q}$, $\ {\widetilde{\!\!\widetilde{{T}}}}$, and $\ {\widetilde{\!\!\widetilde{G}}}$ are given in the Appendix. By inserting $\widehat J$, ${\widetilde J}_{\bm{\omega\omega}}$, ${\:{\J}}_{\bm{\omega\omega}}$, ${\:{\J}}$, $\widetilde{\bm{U}}$, ${\widetilde U}'$, and ${\widehat{\bm{\ell}}}'- {\widetilde{\bm{\ell}}}'$ into (\ref{E.19}) and (\ref{E.4.17}), one obtains the required quantities $\gamma$ and $\rho$ for Barndorff-Nielsen's and Skovgaard's adjustments. Now, one is able to compute the modified statistics $r^*$, $LR^{*}$, and $LR^{**}$.

Computer packages that perform simple operations on matrices and vectors can be used to calculate $\gamma$ and $\rho$. Note that $\gamma$ and $\rho$ depend on the model through ${\bm \mu}_i$, $P_i$, $\Sigma_i$ and $\Sigma_i^{-1}$. The dependence on the specific distribution of ${\bm Y}$ in the class of elliptical distributions occurs through $W_{g}$. 

%------------------------------------- SIMULAÇÃO --------------------------------------------
\section{Simulation study}
\label{sec4}

\inic In this section, we present the results of Monte Carlo simulation experiments in which we evaluate the finite sample performances of the signed likelihood ratio test ($r$) and the standard likelihood ratio test ($LR$) and their corrected versions $r^*$, $LR^*$, and $LR^{**}$. The simulations are based on the univariate nonlinear model and the multivariate mixed linear model when ${\bm Y}_i$ follows a normal distribution, a Student $t$ distribution with $\nu = 3$ degrees of freedom, or a power exponential distribution with shape parameter $\lambda = 0.9$. All simulations are performed using the \texttt{Ox} matrix programming language (Doornik, 2013). The number of Monte Carlo replications is 10,000 (ten thousand). The tests are carried out at the following nominal levels: $\alpha = 1\%, 5\%, 10\%$.

First consider the nonlinear model defined in (\ref{MainModel1}) with
\begin{equation}\label{nonlinear-model}
{\mu}_{i} = \frac{1}{1 + \beta_0 + \beta_1 x_{i1} + \beta_2 x_{i2} + \beta_3 x_{i2}^2}, \quad i = 1,\ldots, n
\end{equation}
(model 1). We test ${\mathcal H}_{0}: \beta_3 \geq 0$ against ${\mathcal H}_{1}: \beta_3 < 0$ and ${\mathcal H}_{0}: (\beta_2, \beta_3)^{\top} = (0, 0)^{\top}$ against ${\mathcal H}_{1}: (\beta_2, \beta_3)^{\top} \neq (0, 0)^{\top}$. The values of the covariates ${x}_{i1}$ and ${x}_{i2}$ are taken as random draws from the standard uniform distribution ${\mathcal U}(0,1)$ and $n = 15, 25, 35, 50$. The parameter values are set at $\beta_0 = 0.5$, $\beta_1 = 0.2$, $\beta_2 = 0$, $\beta_3 = 0$, and $\sigma^2 = 0.005$. For this parameter setting, ${\mu}_i\in (1/(1+0.5+0.2), 1/(1+0.5)) \approx (0.59, 0.67)$ because $x_{ij}\in (0,1)$. This implies that $\sigma^2$ must be very small for the response variable not to be dominated by the random noise. The null rejection rates (say $\widehat p$) of the tests are estimates of the true type I error probabilities with standard error ${\rm se}=\sqrt{\widehat p (1-\widehat p)/10000}$, and 95\% confidence intervals are given by $\widehat p \pm 1.96 \sqrt{\widehat p (1-\widehat p)/10000}$. The null rejection rates and the corresponding standard errors are displayed in Tables $1$ (one-sided tests) and $2$ (two-sided tests) for different sample sizes. 

The simulation results show that the test based on the modified signed likelihood ratio statistic, $r^*$, presents rejection rates closer to the nominal levels than the original version, $r$, in small samples. For instance, for the normal distribution, $n=25$, and $\alpha = 1\%$, the rejection rates are $2.1\%$ $(r)$ and $1.4\%$ $(r^*)$. For $n = 15$, Student $t$ distribution and $\alpha = 5\%$, the rejection rates are  $11.9\%$ $(r)$ and $6.5\%$ $(r^*)$. Additionally, we observe that the test based on the standard likelihood ratio statistic, $LR$, is considerably liberal when the sample size is small, i.e. the rejection rates are much larger than the nominal levels. For instance, for the Student $t$ distribution and $n = 15$, the rejection rates for the $LR$ test equal $5.7\%$ $(\alpha = 1\%)$, $15.7\%$ $(\alpha = 5\%)$, and $25.2\%$ $(\alpha = 10\%)$ (Table 2). The tests based on the modified versions, $LR^{*}$ and $LR^{**}$, present rejection rates much closer to the nominal levels than the original version, $LR$. For example, for the Student $t$ distribution, $n=15$ and $\alpha=1\%$, the rejection rates are $5.7\%$ $(LR)$, $1.1\%$ $(LR^{*})$, and  $0.8\%$ $(LR^{**})$. The above mentioned findings are corroborated by comparing the 95\% confidence intervals of the true probability type I error of the different tests.

\begin{table}[!htp]
\hspace{1.6cm} {\caption{Null rejection rates of the tests of ${\mathcal H}_{0}: \beta_3 \geq 0$ against ${\mathcal H}_{1}: \beta_3 < 0$ (standard errors between parentheses); model 1. Entries are percentages.}} 
\vspace{0.5cm} \centering
\begin{tabular}{ccccccccc} 
\hline
\multicolumn{9}{c}{Normal distribution}\\
\hline
&\multicolumn{2}{c}{$\alpha = 1\%$}&&\multicolumn{2}{c}{$\alpha = 5\%$}&&\multicolumn{2}{c}{$\alpha = 10\%$}
\\\cline{2-3}\cline{5-6}\cline{8-9}
{$n$}& $r$ & $r^*$ && $r$  & $r^*$ && $r$  & $r^*$ \\\hline 
  15 &2.8 (0.2) &1.5 (0.1) && 8.8 (0.3) &6.1 (0.2) &&14.4 (0.4) &11.0 (0.3)  \\
  25 &2.1 (0.1) &1.4 (0.1) && 7.6 (0.3) &6.2 (0.2) &&12.8 (0.3) &11.1 (0.3)  \\ 
  35 &1.7 (0.1) &1.5 (0.1) && 6.5 (0.3) &5.7 (0.2) &&12.0 (0.3) &10.8 (0.3)  \\   
  50 &1.3 (0.1) &1.1 (0.1) && 5.6 (0.2) &5.1 (0.2) &&11.1 (0.3) &10.3 (0.3)  \\
\hline
\multicolumn{9}{c}{Student $t$ distribution}\\
\hline
&\multicolumn{2}{c}{$\alpha = 1\%$}&&\multicolumn{2}{c}{$\alpha = 5\%$}&&\multicolumn{2}{c}{$\alpha = 10\%$}
\\\cline{2-3}\cline{5-6}\cline{8-9}
{$n$}& $r$ & $r^*$ && $r$  & $r^*$ && $r$  & $r^*$ \\\hline 
  15 &4.3 (0.2) &1.8 (0.1) &&11.9 (0.3) &6.5 (0.3) &&18.1 (0.4) &12.2 (0.3)  \\
  25 &2.3 (0.2) &1.7 (0.1) && 7.9 (0.3) &6.2 (0.2) &&13.6 (0.3) &11.5 (0.3)  \\ 
  35 &1.9 (0.1) &1.8 (0.1) && 7.1 (0.3) &6.2 (0.2) &&12.4 (0.3) &11.5 (0.3)  \\   
  50 &1.6 (0.1) &1.4 (0.1) && 6.0 (0.2) &5.4 (0.2) &&11.4 (0.3) &10.6 (0.3)  \\
\hline
\multicolumn{9}{c}{Power exponential distribution}\\
\hline
&\multicolumn{2}{c}{$\alpha = 1\%$}&&\multicolumn{2}{c}{$\alpha = 5\%$}&&\multicolumn{2}{c}{$\alpha = 10\%$}
\\\cline{2-3}\cline{5-6}\cline{8-9}
{$n$}& $r$ & $r^*$ && $r$  & $r^*$ && $r$  & $r^*$ \\\hline 
  15 &2.9 (0.2) &1.9 (0.1) && 9.2 (0.3) &6.5 (0.3) &&15.1 (0.4) &11.9 (0.3)  \\
  25 &1.6 (0.1) &1.4 (0.1) && 7.2 (0.3) &5.8 (0.2) &&12.7 (0.3) &11.5 (0.3)  \\ 
  35 &1.6 (0.1) &1.6 (0.1) && 6.7 (0.3) &6.2 (0.2) &&12.2 (0.3) &11.5 (0.3)  \\   
  50 &1.4 (0.1) &1.2 (0.1) && 6.1 (0.2) &5.7 (0.2) &&11.5 (0.3) &11.0 (0.3)  \\
\hline
\end{tabular}
\label{tab.1}
\end{table}

\begin{table}[!htp]
\hspace{1.6cm} {\caption{Null rejection rates of the tests of ${\mathcal H}_{0}: (\beta_2, \beta_3)^{\top} = (0, 0)^{\top}$ against ${\mathcal H}_{1}: (\beta_2, \beta_3)^{\top} \neq (0, 0)^{\top}$ (standard errors between parentheses); model 1. Entries are percentages.}} 
\vspace{0.5cm} \centering
{
\footnotesize
\begin{tabular}{c|ccc|ccc|ccc} 
\hline
\multicolumn{10}{c}{Normal distribution}\\
\hline
&\multicolumn{3}{c|}{$\alpha = 1\%$}&\multicolumn{3}{c|}{$\alpha = 5\%$}&\multicolumn{3}{c}{$\alpha = 10\%$}
\\\cline{2-4}\cline{5-7}\cline{8-10}
{$n$}& $LR$ & $LR^*$& $LR^{**}$ & $LR$  & $LR^*$& $LR^{**}$ & $LR$  & $LR^*$& $LR^{**}$ \\\hline 
  15 &3.3 (0.2) &1.0 (0.1) &0.9 (0.1) &11.0   (0.3) &5.2 (0.2) &4.9 (0.2) &18.4 (0.4) & 10.5  (0.3) &\ 9.9\ (0.3)  \\
  25 &2.4 (0.2) &1.3 (0.1) &1.2 (0.1) &\ 8.6\ (0.3) &5.9 (0.2) &5.7 (0.2) &15.3 (0.4) & 10.8  (0.3) & 10.6  (0.3)  \\ 
  35 &2.0 (0.1) &1.1 (0.1) &1.1 (0.1) &\ 7.8\ (0.3) &5.6 (0.2) &5.5 (0.2) &13.8 (0.3) & 10.8  (0.3) & 10.7  (0.3)  \\   
  50 &1.3 (0.1) &0.9 (0.1) &0.9 (0.1) &\ 6.2\ (0.2) &4.9 (0.2) &4.9 (0.2) &11.7 (0.3) &\ 9.8\ (0.3) &\ 9.8\ (0.3)  \\
\hline
\multicolumn{10}{c}{Student $t$ distribution}\\
\hline
&\multicolumn{3}{c|}{$\alpha = 1\%$}&\multicolumn{3}{c|}{$\alpha = 5\%$}&\multicolumn{3}{c}{$\alpha = 10\%$}
\\\cline{2-4}\cline{5-7}\cline{8-10}
{$n$}& $LR$ & $LR^*$& $LR^{**}$ & $LR$  & $LR^*$& $LR^{**}$ & $LR$  & $LR^*$& $LR^{**}$ \\\hline 
  15 &5.7 (0.2) &1.1 (0.1) &0.8 (0.1) & 15.7  (0.4) &5.1 (0.2) &4.2 (0.2) &25.2 (0.4) &10.1 (0.3) &\ 8.5\ (0.3) \\
  25 &2.6 (0.2) &1.0 (0.1) &0.9 (0.1) &\ 9.2\ (0.3) &5.2 (0.2) &5.0 (0.2) &16.4 (0.4) &10.3 (0.3) &\ 9.9\ (0.3) \\
  35 &2.0 (0.1) &1.0 (0.1) &1.0 (0.1) &\ 7.9\ (0.3) &5.2 (0.2) &5.1 (0.2) &13.9 (0.4) &10.2 (0.3) & 10.1  (0.3) \\
  50 &1.5 (0.1) &0.9 (0.1) &0.9 (0.1) &\ 6.7\ (0.3) &5.0 (0.2) &4.9 (0.2) &12.8 (0.3) &10.1 (0.3) & 10.0  (0.3) \\
\hline
\multicolumn{10}{c}{Power exponential distribution}\\
\hline
&\multicolumn{3}{c|}{$\alpha = 1\%$}&\multicolumn{3}{c|}{$\alpha = 5\%$}&\multicolumn{3}{c}{$\alpha = 10\%$}
\\\cline{2-4}\cline{5-7}\cline{8-10}
{$n$}& $LR$ & $LR^*$& $LR^{**}$ & $LR$  & $LR^*$& $LR^{**}$ & $LR$  & $LR^*$& $LR^{**}$ \\\hline 
  15 &3.8 (0.2) &1.3 (0.1) &1.2 (0.1) & 12.1  (0.3) &5.8 (0.2) &5.3 (0.2) &19.7 (0.4) & 11.1  (0.3) & 10.4  (0.3) \\
  25 &2.1 (0.1) &1.0 (0.1) &0.9 (0.1) &\ 7.9\ (0.3) &5.0 (0.2) &4.9 (0.2) &14.5 (0.4) &\ 9.9\ (0.3) &\ 9.7\ (0.3) \\
  35 &1.9 (0.1) &1.1 (0.1) &1.1 (0.1) &\ 7.4\ (0.3) &5.5 (0.2) &5.4 (0.2) &14.1 (0.4) & 11.0  (0.3) & 10.8  (0.3) \\
  50 &1.3 (0.1) &0.9 (0.1) &0.9 (0.1) &\ 6.3\ (0.2) &4.9 (0.2) &4.9 (0.2) &11.9 (0.3) & 10.0  (0.3) &\ 9.9\ (0.3) \\
\hline
\end{tabular}
}
\label{tab.2}
\end{table}

Figures 1 and 2 depict curves of relative $p$-values discrepancies {\it versus} the corresponding asymptotic $p$-values for the tests that use $r$ and $r^*$ (Figure 1), and $LR$, $LR^{*}$, and $LR^{**}$ (Figure 2) for $n = 15$ under normal, Student $t$ and power exponential distributions. The relative $p$-value discrepancy is defined by the difference between the exact and the asymptotic $p$-values divided by the asymptotic $p$-value. The closer to zero the curve is the better the asymptotic approximation. The figures clearly suggest that the modified statistics are much better approximated by the respective asymptotic distributions than the unmodified ones.

\begin{figure}[!htp]  
\begin{center}
\includegraphics[height = 6cm, width=10cm]{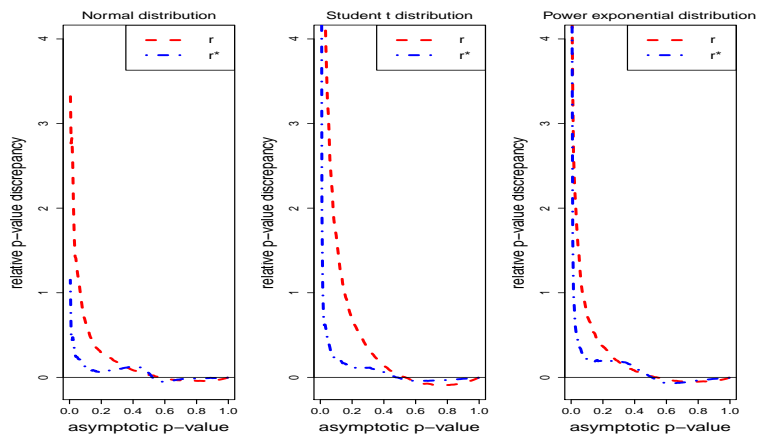}
\caption{Relative $p$-value discrepancy plots for the test of ${\mathcal H}_{0}: \beta_3 \geq 0$ against ${\mathcal H}_{1}: \beta_3 < 0$ for $n=15$; model 1.}
\label{Fig1}
\end{center}
\end{figure}

\begin{figure}[!htp]  
\begin{center}
\includegraphics[height = 6cm, width=10cm]{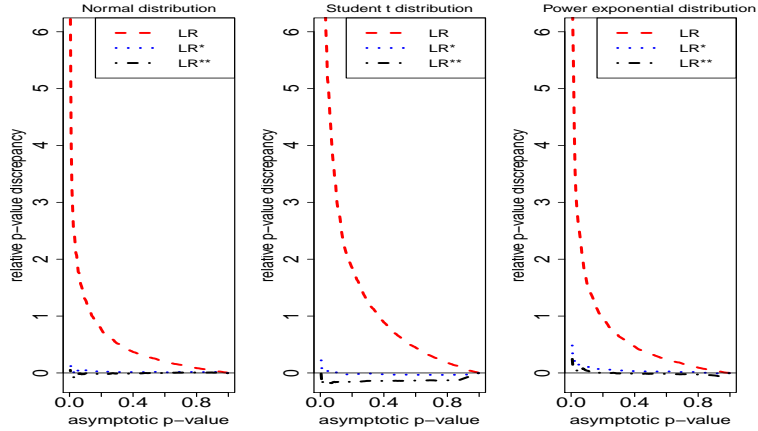}
\caption{Relative $p$-value discrepancy plots for the test of ${\mathcal H}_{0}: (\beta_2, \beta_3)^{\top} = (0, 0)^{\top}$ against ${\mathcal H}_{1}: (\beta_2, \beta_3)^{\top} \neq (0, 0)^{\top}$ for $n=15$; model 1.}
\label{Fig2}
\end{center}
\end{figure}

We conclude that the modified versions of the signed likelihood ratio test and standard likelihood ratio test have better performance than the original tests for small and moderate sample sizes. Although the tests based on the modified statistics, $LR^*$ and $LR^{**}$, present similar results, the test based on $LR^{**}$ had better performance (in the majority of the cases) than the one based on $LR^*$ in our simulations. However, a comparison based on 95\% confidence intervals is not able to distinguish between the two modified tests in most of the cases.

We now consider the mixed linear model
\begin{equation}\label{mixed-model}
\bm{Y}_i =  X_i {\bm{\beta}} + Z_i \bm{b}_i + \bm{e}_i, \quad i = 1,\ldots, n,
\end{equation}
where $\bm{Y}_i$ is the $q_i \times 1$ response vector with $q_i$ randomly chosen from $\{1, 2, 3, 4, 5\}$, $X_i = [{\bm 1} \ \ {\bm x}_{i1} \ \ {\bm x}_{i2} \ \ {\bm x}_{i3}  \ \ {\bm x}_{i4}]$ is a matrix of nonstochastic covariates $(q_i \times 5)$, and ${Z}_i = [{\bm 1} \ \ {\bm x}_{i1}]$ is a matrix of known constants $(q_i \times 2)$ (model 2). The vector ${\bm x}_{i1}$ is composed by the first $q_i$ elements of $(5, 10, 15, 30, 60)^\top$; ${\bm x}_{i2}, {\bm x}_{i3}$, and ${\bm x}_{i4}$ are vectors of dummy variables. 
The vector of the fixed effects parameters is $\bm{\beta} = (\beta_0, \beta_1, \beta_2, \beta_3, \beta_4)^\top$. 
Assume that $(\bm{e}_i, {\bm b}_i)^\top \sim El_{q_i}(\bm{0}, S_i)$, where 
\begin{equation*}\label{MatrizD}
S_i = \left[\begin{array}{cc} \sigma^2 I_{q_i} & 0_{q_i \times 2}  \\ 0_{2 \times q_i} & \Delta({\bm{\gamma}})  \\ \end{array}\right], \ \ \Delta({\bm{\gamma}}) = \left[\begin{array}{cc} \gamma_1 & \gamma_2  \\ \gamma_2 & \gamma_3  \\ \end{array}\right],
\end{equation*}
and ${\bm{\gamma}} = (\gamma_1, \gamma_2, \gamma_3)^\top$. Here, $0_{q_i \times 2}$ is null matrix of dimension $q_i \times 2$.
Therefore, the marginal distribution of the observed vector is $\bm{Y}_i \sim El_{q_i}\left({\bm \mu}_{i}; \Sigma_i\right)$, where ${\bm \mu}_{i} = X_i {\bm{\beta}}$ and ${\Sigma}_i = Z_i \Delta({\bm{\gamma}}) Z_i^{\top} + \sigma^2 I_{q_i}$. Note that model (\ref{mixed-model}) is a special case of (\ref{MainModel}). Here, the vector of unknown parameters is ${\bm{\theta}} = ({\bm{\beta}}^\top, {\bm{\gamma}}^\top, \sigma^2)^\top$. The sample sizes considered are $n = 16, 24, 32, 40$, and $48$. We test ${\mathcal H}_0: {\bm \psi} = {\bm 0}$ against ${\mathcal H}_1: {\bm \psi} \neq {\bm 0}$, where ${\bm \psi} = (\beta_2, \beta_3, \beta_4)^\top$. The parameter values are $\beta_0 = 0.7$, $\beta_1 = 0.5$, $\beta_2 = \beta_3 = \beta_4 = 0$, $\gamma_1 = 500$, $\gamma_2 = 2$, $\gamma_3 = 200$, and $\sigma^2 = 5$. 

The null rejection rates of the tests are displayed in Table 3. We note that the likelihood ratio test is liberal. 
For instance, when ${\bm Y}_i$ follows a Student $t$ distribution, $n = 16$, and $\alpha = 10\%$, its rejection rate exceeds $26\%$. 
The tests based on the modified statistic, $LR^{*}$ and $LR^{**}$, present rejection rates closer to the nominal levels than the original version, $LR$. For example, the rejection rates when ${\bm Y}_i$ follows a Student $t$ distribution, $n = 16$, and $\alpha = 1\%$, are $6.3\%$ ($LR$), $1.1\%$ ($LR^{*}$), and $0.9\%$ ($LR^{**}$). 

Figure 3 presents curves of relative $p$-values discrepancies {\it versus} the corresponding asymptotic $p$-values for the statistics $LR$, $LR^{*}$, and $LR^{**}$ for $n = 16$ under the distributions considered. It is evident that the modified statistics are much better approximaded by the reference distributions than the original statistics.

In conclusion, the simulations suggest that the modified versions of the standard and signed likelihood ratio tests perform better than the original tests for small and moderate sample sizes. 

\begin{table}[!htp]
\hspace{1.6cm} {\caption{Null rejection rates (given in percentages) of the tests of ${\mathcal H}_0: {\bm \psi} = {\bm 0}$ against ${\mathcal H}_1: {\bm \psi} \neq {\bm 0}$, where ${\bm \psi} = (\beta_2, \beta_3, \beta_4)^\top$ (standard errors between parentheses); model 2. Entries are percentages.}} 
\vspace{0.5cm} \centering
{
\footnotesize
\begin{tabular}{c|ccc|ccc|ccc} 
\hline
\multicolumn{10}{c}{Normal distribution}\\
\hline
&\multicolumn{3}{c|}{$\alpha = 1\%$}&\multicolumn{3}{c|}{$\alpha = 5\%$}&\multicolumn{3}{c}{$\alpha = 10\%$}
\\\cline{2-4}\cline{5-7}\cline{8-10}
{$n$}& $LR$ & $LR^*$& $LR^{**}$ & $LR$  & $LR^*$& $LR^{**}$ & $LR$  & $LR^*$& $LR^{**}$ \\\hline 
  16 &5.8 (0.2) &1.3 (0.1) &1.0 (0.1) & 16.1  (0.4) &5.8 (0.2) &4.9 (0.2) &24.6 (0.4) &10.9 (0.3) &\ 9.8\ (0.3) \\
  24 &3.2 (0.2) &1.3 (0.1) &1.2 (0.1) & 10.8  (0.3) &5.3 (0.2) &5.1 (0.2) &18.2 (0.4) &10.6 (0.3) & 10.2  (0.3) \\ 
  32 &2.5 (0.2) &1.1 (0.1) &1.1 (0.1) &\ 9.2\ (0.3) &5.4 (0.2) &5.2 (0.2) &15.7 (0.4) &10.5 (0.3) & 10.3  (0.3) \\   
  40 &2.0 (0.1) &1.0 (0.1) &1.0 (0.1) &\ 8.0\ (0.3) &5.0 (0.2) &4.9 (0.2) &14.7 (0.4) &10.3 (0.3) & 10.1  (0.3) \\ 
  48 &1.9 (0.1) &1.0 (0.1) &1.0 (0.1) &\ 7.3\ (0.3) &4.9 (0.2) &4.9 (0.2) &13.7 (0.3) &10.2 (0.3) & 10.1  (0.3) \\ 
\hline
\multicolumn{10}{c}{Student $t$ distribution}\\
\hline
&\multicolumn{3}{c|}{$\alpha = 1\%$}&\multicolumn{3}{c|}{$\alpha = 5\%$}&\multicolumn{3}{c}{$\alpha = 10\%$}
\\\cline{2-4}\cline{5-7}\cline{8-10}
{$n$}& $LR$ & $LR^*$& $LR^{**}$ & $LR$ & $LR^*$& $LR^{**}$ & $LR$  & $LR^*$& $LR^{**}$ \\\hline 
  16 &6.3 (0.2) &1.1 (0.1) &0.9 (0.1) & 17.4  (0.4) &5.4 (0.2) &4.6 (0.2) &26.8 (0.4) & 10.8  (0.3) &\ 9.1\ (0.3) \\ 
  24 &3.5 (0.2) &1.0 (0.1) &0.9 (0.1) & 11.4  (0.3) &5.4 (0.2) &5.1 (0.2) &19.1 (0.4) & 10.6  (0.3) & 10.1  (0.3) \\  
  32 &2.7 (0.2) &1.1 (0.1) &1.0 (0.1) &\ 9.3\ (0.3) &5.3 (0.2) &5.0 (0.2) &16.0 (0.4) & 10.1  (0.3) &\ 9.9\ (0.3) \\    
  40 &2.2 (0.2) &1.1 (0.1) &1.1 (0.1) &\ 8.4\ (0.3) &5.4 (0.2) &5.2 (0.2) &14.8 (0.4) & 10.3  (0.3) & 10.1  (0.3) \\ 
  48 &1.8 (0.1) &1.0 (0.1) &0.9 (0.1) &\ 7.2\ (0.3) &4.6 (0.2) &4.5 (0.2) &13.4 (0.3) &\ 9.9\ (0.3) &\ 9.8\ (0.3) \\ 
\hline
\multicolumn{10}{c}{Power exponential distribution}\\
\hline
&\multicolumn{3}{c|}{$\alpha = 1\%$}&\multicolumn{3}{c|}{$\alpha = 5\%$}&\multicolumn{3}{c}{$\alpha = 10\%$}
\\\cline{2-4}\cline{5-7}\cline{8-10}
{$n$}& $LR$ & $LR^*$& $LR^{**}$ & $LR$ & $LR^*$& $LR^{**}$ & $LR$  & $LR^*$& $LR^{**}$ \\\hline 
  16 &6.0 (0.2) &1.2 (0.1) &1.0 (0.1) & 16.2  (0.4) &6.2 (0.2) &5.3 (0.2) &24.8 (0.4) &11.8 (0.3) &10.5 (0.3) \\
  24 &3.2 (0.2) &1.1 (0.1) &1.1 (0.1) & 10.9  (0.3) &5.3 (0.2) &5.1 (0.2) &18.1 (0.4) &10.6 (0.3) &10.2 (0.3) \\  
  32 &2.3 (0.2) &1.1 (0.1) &1.0 (0.1) &\ 9.4\ (0.3) &5.1 (0.2) &5.0 (0.2) &16.0 (0.4) &10.9 (0.3) &10.6 (0.3) \\  
  40 &2.0 (0.1) &1.1 (0.1) &1.1 (0.1) &\ 7.9\ (0.3) &5.2 (0.2) &5.1 (0.2) &14.5 (0.4) &10.2 (0.3) &10.0 (0.3) \\ 
  48 &1.9 (0.1) &1.1 (0.1) &1.1 (0.1) &\ 7.4\ (0.3) &5.4 (0.2) &5.4 (0.2) &13.7 (0.3) &10.4 (0.3) &10.3 (0.3) \\ 
\hline
\end{tabular}
}
\label{tab.3}
\end{table}

\begin{figure}[!htp]  
\begin{center}
\includegraphics[height = 6cm, width=10cm]{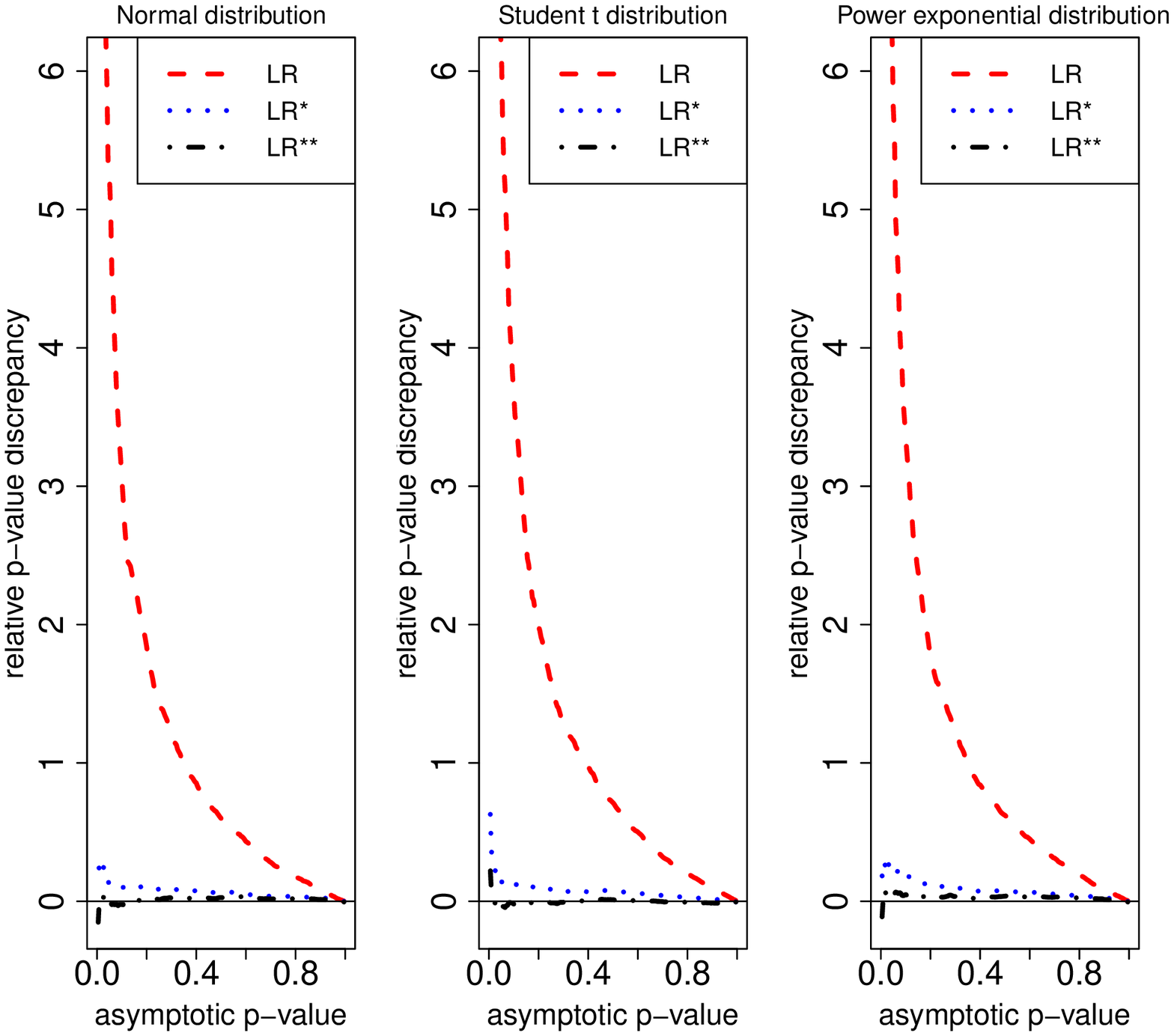}
\caption{Relative $p$-value discrepancy plots for the test of ${\mathcal H}_0: {\bm \psi} = {\bm 0}$ against ${\mathcal H}_1: {\bm \psi} \neq {\bm 0}$, where ${\bm \psi} = (\beta_2, \beta_3, \beta_4)^\top$ for $n=16$; model 2.}
\label{Fig3}
\end{center}
\end{figure}

%------------------------------------- APLICAÇÃO --------------------------------------------
\section{Applications}
\label{sec5}

\subsection{Fluorescent lamp data}\label{lamp}

\inic The data used in this section were presented by Rosillo et al. (2010, Table 5). The authors analyzed the lifetime of $n=14$ fluorescent lamps in photovoltaic systems using an analytical model whose goal is to assist in improving ballast design and extending the lifetime of fluorescent lamps.

We consider the nonlinear model (\ref{nonlinear-model}) where the response variable is the observed lifetime/advertised lifetime ($Y$), the covariates correspond to a measure of gas discharge ($x_1$) and the observed voltage/advertised voltage (measure of performance of lamp and ballast -- $x_2$); see Rosillo et al. (2010, eq(15)). The errors are assumed to follow a Student $t$ distribution with $\nu = 4$ degrees of freedom. This model provides a suitable fit for the data (Melo et al, 2015). The maximum likelihood estimates of the parameters (standard errors are given in parentheses) are $\widehat\beta_0 = 33.519 \ (5.082)$, $\widehat\beta_1 = 9.592 \ (4.417)$, $\widehat\beta_2 = -63.501 \ (9.789)$, $\widehat\beta_3 = 29.777 \ (4.710)$, and $\widehat\sigma^2 = 0.006 \ (0.003)$.

The signed likelihood ratio statistic for the one-sided test of ${\mathcal H}_0: \beta_1 \leq 0$ against ${\mathcal H}_1: \beta_1 > 0$ is $r = 2.374$ ($p$-value: $0.009$), and the corrected statistic is $r^{*} = 1.878$ ($p$-value: $0.030$). The unmodified test rejects the null hypothesis at the 1\% significance level, unlike the modified test. In addition, the standard likelihood ratio statistic for testing the two-sided test of
${\mathcal H}_0: \beta_1 = 0$ against ${\mathcal H}_1: \beta_1 \neq 0$ is $LR = 5.634$ ($p$-value: $0.018$) whereas the corrected statistics are $LR^{*} = 3.528$ ($p$-value: $0.060$) and $LR^{**} = 3.282$ ($p$-value: $0.070$). Note the considerable increase in the $p$-values when the 
modified statistics are employed. At the 5\% significance level, the null hypothesis is rejected by the tests based on the modified statistics $LR^{*}$ and $LR^{**}$ but not by the standard likelihood ratio test.   
  
\subsection{Blood pressure data}\label{blood}

\inic We consider a randomly selected sub-sample of the data presented by Crepeau et al.~(1985). Heart attacks were induced in rats exposed to four different low concentrations of halothane; group 1: 0\% (control), group 2: 0.25\%, group 3: 0.50\% and group 4: 1.0\%. Our sample consists of 35 rats. The blood pressure of each rat (in mm Hg) is recorded over different points in time, from 1 to 9 recordings, after the induced heart attack. The goal is to investigate the effect of halothane on the blood pressure.

We consider the mixed linear model (\ref{mixed-model}), where the $j$th element of the response variable ${\bm Y}_i$ is the blood pressure of the $i$th rat at time $j$ with $i = 1, 2, \ldots, n$, $j = 1, 2, \ldots, q_i$ and $q_i \in \{1, 2, 3, 4, 5, 6, 7, 8, 9\}$. The values of ${\bm x}_{i1}$ are obtained of the vector from points in time (in minutes) in which the $i$th rat blood pressure was recorded. This vector is given by $(5, 10, 15, 30, 60, 120, 180, 240)^\top$. Furthermore, ${\bm x}_{i2}$ is a dummy variable that equals 1 if the $i$th rat belongs to group 2 and 0 otherwise. Also, ${\bm x}_{i3}$ and ${\bm x}_{i4}$ equal 1 for groups 3 and 4, respectively. As in Crepeau et al.~(1985), we assume a normal distribution for ${\bm Y}_i$. The hypothesis ${\mathcal H}_0: {\bm \psi} = {\bm 0}$ is to be tested against  ${\mathcal H}_1: {\bm \psi} \neq {\bm 0}$, where ${\bm \psi} = (\beta_2, \beta_3, \beta_4)^\top$. 

The maximum likelihood estimates of the parameters (standard errors are given in parentheses) are $\widehat\beta_0 = 100.60 \ (6.74)$, $\widehat\beta_1 = 0.011 \ (0.012)$, $\widehat\beta_2 = 4.493 \ (9.333)$, $\widehat\beta_3 = -11.032 \ (9.081)$, $\widehat\beta_4 = -23.022 \ (8.913)$, $\widehat\sigma^2 = 97.761 \ (10.660)$, $\widehat\gamma_1 = 483.550 \ (122.770)$, $\widehat\gamma_2 = -0.700 \ (0.297)$, and $\widehat\gamma_3 = 0.002 \ (0.001)$. The likelihood ratio statistic and its modified versions for testing ${\cal H}_0: {\bm \psi} = {\bm 0}$ are: $LR = 7.954$ ($p$-value: $0.047$), $LR^{*} = 6.883$ ($p$-value: $0.075$), and $LR^{**} = 6.844$ ($p$-value: $0.077$). We notice that the null hypothesis is rejected at the 5\% nominal level when one uses the original statistic, but ${\cal H}_0$ is not rejected when the modified likelihood ratio tests are employed. That is, when using the modified tests (at the 5\% nominal level) one concludes that there is not enough evidence that the blood pressure is affected by the administration of halothane at the concentrations considered in the experiment. This conclusion agrees with the analysis of Crepeau et al.~(1985, p.510) and the test of ${\cal H}_0$ based on the full sample with $43$ rats: $LR = 7.162$ ($p$-value: $0.067$), $LR^{*} = 6.613$ ($p$-value: $0.085$), and $LR^{**} = 6.602$ ($p$-value: $0.086$).

%------------------------------------- CONCLUSÃO --------------------------------------------
\section{Concluding remarks}
\label{sec6}

\inic We studied the issue of testing two-sided and one-sided hypotheses in a general multivariate elliptical model.
Some special cases of this model are errors-in-variables models, nonlinear mixed-effects models, heteroscedastic
nonlinear models, among others. In any of these models, the vector of the errors may have any multivariate elliptical distribution. 
In small and moderate-sized samples the distributions of the standard and signed likelihood ratio statistics may be far from the
respective reference distributions. As a consequence, the tests may be considerably liberal. We derived modified versions of these 
statistics. Our simulations suggest that the modified statistics closely follow the reference distributions in finite samples. 
The modifications obtained in this paper attenuate the liberal behavior of the original tests.

\section*{Acknowledgement}
We gratefully acknowledge the financial support from CNPq and FAPESP. 
The authors thank the anonymous referee and the Associate Editor for helpful comments and suggestions.

%------------------------------------- APÊNDICE --------------------------------------------
\begin{appendix}
\section*{Appendix. The observed information matrix and derivatives with respect to the data}
\label{secA}

\inic The $(r,s)$th element of the observed information matrix $J$ is given by 
$$\sum_{i=1}^{n}\Big\{\bm{T}_{i(r)}^{\top} \Sigma_i^{-1} \bm{d}_{i(s)} + \tr(B_{i(r)} A_{i(s)}) + E_{i(rs)}\Big\},$$ where
\begin{equation*}
\begin{split}
\bm{T}_{i(r)}^{\top} &= - \dot{v}_i \left(\bm{z}_i^{\top} A_{i(r)} \bm{z}_i\right) \bm{z}_i^{\top} + 2 \dot{v}_i \left(\bm{d}_{i(r)}^{\top} \Sigma_i^{-1} \bm{z}_i\right) \bm{z}_i^{\top} + v_i \bm{d}_{i(r)}^{\top} + v_i \bm{z}_i^{\top} \Sigma_i^{-1} C_{i(r)}, \\
B_{i(r)}        &= - \dot{v}_i \left(\bm{d}_{i(r)}^{\top} \Sigma_i^{-1} \bm{z}_i\right) \bm{z}_i \bm{z}_i^{\top} + \frac{1}{2}\dot{v}_i \left(\bm{z}_i^{\top} A_{i(r)} \bm{z}_i\right) \bm{z}_i \bm{z}_i^{\top} - v_i \bm{z}_i \bm{d}_{i(r)}^{\top} - \frac{1}{2} C_{i(r)},\\
E_{i(rs)}       &= - \frac{1}{2} \tr\left[A_{i(sr)} \left(\Sigma_i - v_i \bm{z}_i \bm{z}_i^{\top}\right)\right] - v_i \bm{z}_i^{\top} \Sigma_i^{-1} \bm{d}_{i(sr)},\\
\end{split}
\end{equation*}
with $\bm{z}_i = {\bm Y}_i - \bm{\mu}_i = \widehat{P}_i \bm{a}_i + \widehat{\bm{\mu}}_i - \bm{\mu}_i$, $\dot{v}_i = - 2 W_g'(u_i)$, $W_g'(u) = d W_g(u)/du$, and $A_{i(sr)} = \partial A_{i(s)}/\partial\theta_r = - 2 A_{i(r)} C_{i(s)} \Sigma_i^{-1} - \Sigma_i^{-1} (\partial C_{i(s)}/\partial \theta_r)\Sigma_i^{-1}$. For the $N_{q_i}({{\bm \mu}_i}, {\Sigma_i})$ distribution, $v_i = 1$ and $\dot{v}_i = 0$. For the $t_{q_i}({{\bm \mu}_i}, {\Sigma_i}, \nu)$ distribution, we have $v_i = (\nu + q_i)/(\nu + u_i)$ and $\dot{v}_i = -(\nu + q_i)/(\nu + u_i)^2$. Additionally, for the $PE_{q_i}({{\bm \mu}_i}, {\Sigma_i}, \lambda)$ distribution, we have $v_i = \lambda u_i^{\lambda - 1}$ and $\dot{v}_i = \lambda (\lambda - 1) u_i^{\lambda - 2}$. 

Hence, the observed information matrix can be written as in (\ref{ScoreFisher2}), where 
$$G_{rs} = \sum_{i=1}^{n}\left[\tr(B_{i(r)} A_{i(s)}) + E_{i(rs)}\right].$$

We now turn to the derivatives with respect to the sample space required to compute Barndorff-Nielsen's and Skovgaard's adjustments.
% are $\ell'$ and $U'$. 
From (\ref{E.ell}), the $r$th element of the vector ${\bm{\ell}}'$ is 
\begin{equation*}\label{C.2}
\begin{split}
\ell_r' &= \sum_{i=1}^{n} \left(\bm{a}_i^{\top} \widehat{P}_{i(r)}^{\top} + \widehat{\bm{d}}_{i(r)}^{\top}\right) \Sigma_i^{-1} \left(-v_i \bm{z}_i\right),
\end{split}
\end{equation*}
where $\widehat{\bm{d}}_{i(r)} = \partial {\widehat{\bm{\mu}}}_{i}/\partial{\widehat{\theta}}_r$ and ${\widehat P}_{i(r)} = \partial {\widehat P}_{i}/\partial{\widehat{\theta}}_r$. The derivatives of $P_i$ with respect to the parameters may be obtained by using the algorithm proposed by Smith (1995). 	  

The $(r,s)$th element of the matrix $U'$ is
$$%\label{C.5}
U_{rs}' = \sum_{i=1}^{n} \left(\bm{a}_i^{\top} \widehat{P}_{i(s)}^{\top} + \widehat{\bm{d}}_{i(s)}^{\top}\right) \Sigma_i^{-1} \Big[2 \dot{v}_i \bm{z}_i \bm{d}_{i(r)}^{\top} \Sigma_i^{-1} \bm{z}_i + v_i \bm{d}_{i(r)} - \dot{v}_i \bm{z}_i \bm{z}_i^{\top} A_{i(r)} \bm{z}_i + v_i C_{i(r)} \Sigma_i^{-1} \bm{z}_i\Big].
$$
In matrix notation ${\bm{\ell}}'$ and $U'$ can be written as 
\begin{equation*}
{\bm{\ell}}' = \widehat{R}^\top \Sigma^{-1} {\bm z}^*, \ \ U' = Q^\top \Sigma^{-1} \widehat{R},
\end{equation*}
where $\widehat{R}^\top = (\widehat{R}_1^{\top}, \ldots, \widehat{R}_n^{\top})$, ${{\bm z}}^* = (- {v}_1 {{\bm z}}_1^{{\top}}, \ldots, - {v}_n {{\bm z}}_n^{{\top}})^{\top}$, $Q = (Q_1^{\top}, \ldots, Q_n^{\top})^\top$, $\widehat{R}_i = (\widehat{\bm R}_{i(1)}, \ldots, \break \widehat{\bm R}_{i(p)})$, and $Q_i = (\bm{Q}_{i(1)}, \ldots, \bm{Q}_{i(p)})$, with $\widehat{\bm{R}}_{i(r)} = \widehat{P}_{i(r)} {\bm a}_i + \widehat{\bm{d}}_{i(r)}$ and
$$
\bm{Q}_{i(r)}           = 2 \dot{v}_i \bm{z}_i \bm{d}_{i(r)}^{\top} \Sigma_i^{-1} \bm{z}_i + v_i \bm{d}_{i(r)} - \dot{v}_i \bm{z}_i \bm{z}_i^{\top} A_{i(r)} \bm{z}_i + v_i C_{i(r)} \Sigma_i^{-1} \bm{z}_i. \\
$$
Therefore, $\widehat{\bm{\ell}}'$, $\widetilde{\bm{\ell}}'$ and $\widetilde{U}'$ can be written as in (\ref{Deriv}), where $\widetilde{v}_i = -2 W_{g}(\widetilde{u}_{i})$, $\widetilde{\dot{v}}_i = -2 W_{g}'(\widetilde{u}_{i})$, $\widehat{{\bm z}}_i = {\widehat P}_i {\bm a}_i$, $\widetilde{{\bm z}}_i = {\widehat P}_i {\bm a}_i + {\widehat{\bm \mu}}_i - \widetilde{{\bm \mu}}_i$, and $\widetilde{u}_i = \left({\widehat P}_i {\bm a}_i + {\widehat{\bm \mu}}_i - \widetilde{{\bm \mu}}_i\right)^{\top} \widetilde{\Sigma}_i^{-1} \left({\widehat P}_i {\bm a}_i + {\widehat{\bm \mu}}_i - \widetilde{{\bm \mu}}_i\right)$. 

Finally, the matrix $\:\:{\J}$ required to compute Skovgaard's adjustment is defined in (\ref{Deriv}), where
 \ \ ${\widetilde{\!\!\widetilde{T}}} = \left(\ \ {\widetilde{\!\!\widetilde{T}}}_1^\top, \ldots, \ \ {\widetilde{\!\!\widetilde{T}}}_n^\top\right)^{\top}$, $\ {\widetilde{\!\!\widetilde{T}}}_i = \left(\ \ {\widetilde{\!\!\widetilde{\bm T}}}_{i(1)}, \ldots, \ \ {\widetilde{\!\!\widetilde{\bm T}}}_{i(p)}\right)$, and the $(r,s)$th element of $\ {\widetilde{\!\!\widetilde{G}}}$ is given by 
$${\widetilde{\!\!\widetilde{G}}}_{rs} = \sum_{i=1}^{n}\left[\tr\left(\ {\widetilde{\!\!\widetilde{B}}}_{i(r)} \widetilde{A}_{i(s)}\right) + \ {\widetilde{\!\!\widetilde{E}}}_{i(rs)}\right]\:,$$ with
\begin{equation*}
\begin{split}
\ {\widetilde{\!\!\widetilde{\bm T}}}_{i(r)} &= - \ {\widetilde{\!\!\widetilde{\dot{v}}}}_{i} \left({\bm a}_i^{\top} {\widetilde P}_i^{\top} \widetilde{A}_{i(r)} {\widetilde P}_i {\bm a}_i\right) {\bm a}_i^{\top} {\widetilde P}_i^{\top} + 2 \ {\widetilde{\!\!\widetilde{\dot{v}}}}_{i} \left(\widetilde{\bm{d}}_{i(r)}^{\top} \widetilde{\Sigma}_i^{-1} {\widetilde P}_i {\bm a}_i \right) {\bm a}_i^{\top} {\widetilde P}_i^{\top} + \ {\widetilde{\!\!\widetilde{v}}}_{i} \widetilde{\bm{d}}_{i(r)}^{\top} + \ {\widetilde{\!\!\widetilde{v}}}_{i} {\bm a}_i^{\top} {\widetilde P}_i^{\top} \widetilde{\Sigma}_i^{-1} \widetilde{C}_{i(r)},\\ \\
\ {\widetilde{\!\!\widetilde{B}}}_{i(r)} &= - \ {\widetilde{\!\!\widetilde{\dot{v}}}}_{i} \left(\widetilde{\bm{d}}_{i(r)}^{\top} \widetilde{\Sigma}_i^{-1} {\widetilde P}_i {\bm a}_i \right) {\widetilde P}_i {\bm a}_i {\bm a}_i^{\top} {\widetilde P}_i^{\top} + \frac{1}{2} \ {\widetilde{\!\!\widetilde{\dot{v}}}}_{i} \left({\bm a}_i^{\top} {\widetilde P}_i^{\top} \widetilde{A}_{i(r)} {\widetilde P}_i {\bm a}_i\right) {\widetilde P}_i {\bm a}_i {\bm a}_i^{\top} {\widetilde P}_i^{\top}  - \ {\widetilde{\!\!\widetilde{v}}}_{i} {\widetilde P}_i {\bm a}_i \widetilde{\bm{d}}_{i(r)}^{\top} \\ & \ \ \ - \frac{1}{2} \widetilde{C}_{i(r)}, \\ \\
\ {\widetilde{\!\!\widetilde{E}}}_{i(rs)}   &= - \frac{1}{2} \tr\left[\widetilde{A}_{i(sr)} \left(\widetilde{\Sigma}_i - \ {\widetilde{\!\!\widetilde{v}}}_{i} {\widetilde P}_i {\bm a}_i {\bm a}_i^{\top} {\widetilde P}_i^{\top}\right)\right] - \ {\widetilde{\!\!\widetilde{v}}}_{i} {\bm a}_i^{\top} {\widetilde P}_i^{\top} \widetilde{\Sigma}_i^{-1} \widetilde{\bm{d}}_{i(sr)}, \\ 
\end{split}
\end{equation*}
where 
$${\widetilde{\!\!\widetilde{v}}}_{i} = -2 W_{g}(\ {\widetilde{\!\!\widetilde{u}}}_{i}), \ \  {\widetilde{\!\!\widetilde{\dot{v}}}}_{i} = -2 W_{g}'(\ {\widetilde{\!\!\widetilde{u}}}_{i}), \ {\rm and} \ \ {\widetilde{\!\!\widetilde{u}}}_{i} = {\bm a}_i^{\top} {\widetilde P}_i^{\top} \widetilde{\Sigma}_i^{-1} {\widetilde P}_i {\bm a}_i.$$ 

\end{appendix}

%------------------------------------- BIBLIOGRAFIA --------------------------------------------

\end{document}